# Smart Vehicle to Grid Interface Project: Electromobility Management System Architecture and Field Test Results

Letterio Zuccaro[*], Alessandro Di Giorgio[*], Francesco Liberati[*], Silvia Canale[*], Andrea Lanna[*], Victor Fernandez Pallares[†], Alejandro Martinez Blanco[†], Raúl Urbano Escobar[†], Jure Ratej[‡], Borut Mehle[‡], Ursula Krisper[§]

[*] DIAG Antonio Ruberti
Sapienza University of Rome, Rome, Italy
Email: {zuccaro, digiorgio, liberati, canale, lanna}@diag.uniroma1.it

[†] CIT Developments, Valencia, Spain
Email: {vfernandez, amartinez, rurbano}@citdev.com

[‡] ETREL, Ljubljana, Slovenia
Email: {jure.ratej, borut.mehle}@etrel.si

[§] Elektro Ljubljana, Ljubljana, Slovenia
Email: ursula.krisper@elektro-ljubljana.si

*Abstract*—This paper presents and discusses the electromobility management system developed in the context of the "SMARTV2G" project, enabling the automatic control of plug-in electric vehicles' (PEVs') charging processes. The paper describes the architecture and the software/hardware components of the electromobility management system. The focus is put in particular on the implementation of a centralized demand side management control algorithm, which allows remote real time control of the charging stations in the field, according to preferences and constraints expressed by all the actors involved (in particular the distribution system operator and the PEV users). The results of the field tests are reported and discussed, highlighting critical issues raised from the field experience.

## I. INTRODUCTION

In a context of obliged continuous optimization of the energy consumption rates in developed societies [1], embedded systems and solutions can play a significant role in the process of transition towards sustainable urban life in European countries. One of the main and most promising technological areas that are expected to be able to contribute in a most relevant way to that overall target is the one constituted by the plug-in electric vehicles (PEVs) [2], [3]. However, previous and current initiatives aiming at a deeper deployment of this environment respectful alternative transport option have had to face serious technological and logistic handicaps that have up to now entailed a noteworthy hurdle for a generalized penetration of this kind of motion technologies.

As technological barriers related to vehicles autonomy seem to be progressively overcome by gradually improving batteries and electric drives, the main technological obstacles for a large scale deployment of PEVs remain the lack of an optimized network of electric vehicle supply equipment (EVSE) – i.e. all accessories, devices, power outlets or apparatuses, including the charging stations (CSs), installed for the purpose of delivering energy to the PEV and allowing communication–, standardized user interfaces and proper charging stations control centers mitigating the effect of simultaneous charging loads on the electricity distribution grid.

As a matter of fact, the main objective targeted by the Smart Vehicle to Grid Interface (SMARTV2G) FP7 European research project [4] aims at connecting the PEVs to the grid by enabling controlled flow of energy and power through safe, secure, energy efficient and convenient transfer of electricity and data. In order to meet the mentioned objective, a control framework enabling PEV charging load control through a demand side management control strategy has been initially proposed in [5], [6] and then refined in [7] and [8][1]. This paper presents an overview over the designed PEV charging control system, in terms of used equipment and control methodology, and then focuses on the achieved field test results, also discussing open points raised during the field experience.

With the term "smart charging", we refer in the following to the possibility of providing the charging service to the PEV users according to their expressed preferences, and in a controlled way, which is suitable and convenient also for the other actors (the charging service provider, the EVSE operator,

---

[1] A different approach analysed in the project was to let the vehicles autonomously decide their charging strategy, see [9].

This work is supported in part by the European Union FP7-2011-ICT-GC SMARTV2G project, under Grant agreement no. 284953.

the distribution system operator (DSO), etc.) involved in the charging process. User preferences are expressed by the user when asking for the charging services and are given by the couple: (i) desired final state of charge (SoC) or, alternatively, required energy to recharge, expressed in kWh and; (ii) available time for charging (i.e. the time within which the final SoC has to be achieved or the required energy has to be delivered).

On the other hand, grid actors, and the DSO in particular, find it convenient that the charging processes (i.e. the power drown by the charging stations) are properly controlled, because to do so would mean for them having the flexibility of employing electromobility controllable load for network-support purposes, balancing of fluctuating renewable energy sources, etc. (thus generating value for themselves and for the other cooperating actors involved – i.e. the PEV users and the operator of the EVSEs).

The remainder of the paper is organized as follows. Section II presents a general introduction to the SMARTV2G system architecture. Section III gives details on the developed SMARTV2G control center. Section IV gives details on the charging station developed in the SMARTV2G project and employed in the field tests. Section V recalls the demand side management algorithm developed for the automatic control of PEV charging loads. Section VI discussed the results of the field tests and, finally, in Section VII the conclusions of the work are given.

## II. SYSTEM ARCHITECTURE

The SMARTV2G project has designed [4] and implemented an electromobility management system to allow actual implementation of smart charging, meaning that:

- An information and communication technology (ICT) - based control center (named "SMARTV2G control center") has been developed to monitor and control the CSs deployed on the field. The PEV users are provided with radio-frequency identification (RFID) cards and the control center takes care of authorization/termination of the charging process.

- The SMARTV2G control center can periodically (each five minutes in the field tests) acquire measurements from the field (i.e. status of the charging stations, energy and power measured by the charging stations currently providing the charging service). As explained next, this feedback of information is crucial to ensure that the load profiles commanded to the charging stations can be properly updated during time.

- An ICT interface (smartphone app or web page) allows PEV users to communicate their willingness to take part in the smart charging program, their user preferences related to the charging session to be started (desired final SoC and time flexibility) and other information to be known by the control center in order to properly control the charging session (the actual/initial SoC and the battery

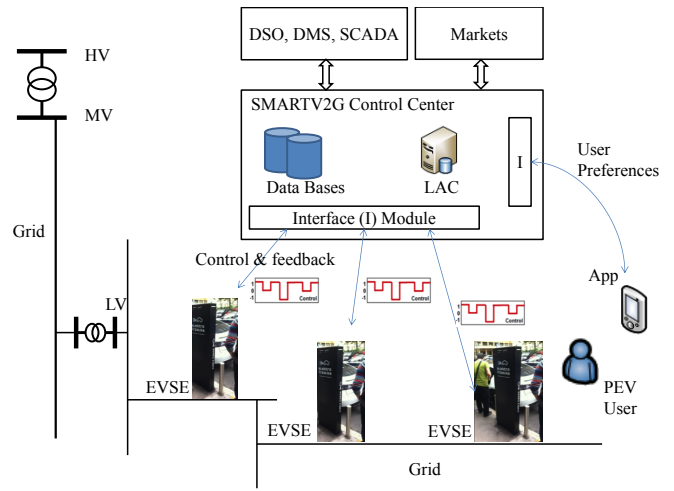

Fig. 1. SMARTV2G reference scenario for smart charging

capacity, which can be derived from vehicle type, to be reported also by the EV user).

- Within the SMARTV2G control center, a software module, named load area controller (LAC), is in charge of solving a demand side management control problem aimed at computing and periodically updating (according to the feedback received from the field) the load profiles to be applied by the charging stations. Such load profiles have always to be compliant, in particular, with the user preferences (as far as grid constraints allows it).

- Smart alternating current CSs developed by ETREL in SMARTV2G project are installed in the field. They are connected with the SMARTV2G control center (via GPRS or Ethernet communication) and support single-phase and tree-phase charging at different power levels. They can be remotely controlled from the control center for various purposes and, in particular for what concerns this paper, in the sense that, according to standard IEC61851 [10], the CS can communicate to the PEV a set point for the value of the current the PEV is allowed to draw in a certain time interval, thus making possible to automatically controlling the EV charging load profiles with proper demand side management algorithms.

Load profiles are computed by the LAC under the two fundamental requirements of: (i) minimizing the costs for charging, even in a dynamic pricing scenario in which the tariff is not flat and may be updated intraday, and (ii) control the aggregated PEV charging load (i.e. the overall power adsorbed by the charging stations) to follow a power reference profile set by a grid actor willing to exploit the flexibility offered by PEVs charging control. Such power reference may be for example set by the DSO on a day-ahead basis, or even updated intraday to face medium/short-term grid needs.

The general reference scenario considered by SMARTV2G is shown in Fig. 1, including also a high-level architecture of the electromobility management system developed in SMARTV2G project. The SMARTV2G control center is discussed in the next section.

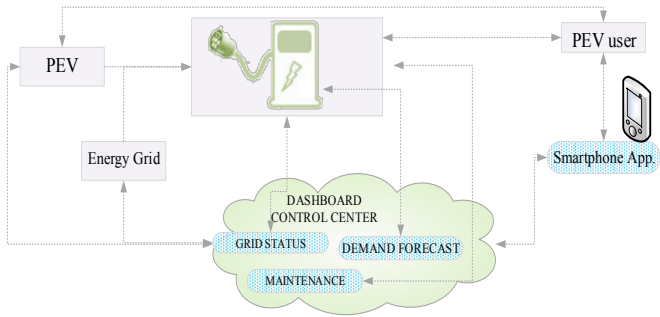

Fig. 2. The SMARTV2G control centre architecture, based on a cloud control center

## III. SMARTV2G CONTROL CENTRE

The architecture of the SMARTV2G control center is based on a cloud solution, which enables the interaction between PEVs, CSs, PEV users and other grid actors (Fig. 2). Several interfaces have been implemented in order to allow the data flowing between these actors. The implemented DashBoard-control center enables operators (such as energy suppliers and business operators) to monitor in real time and also forecast the behavior of the electro mobility grid. Moreover, several proactive functionalities are enabled leveraging the communication protocols established between smart CSs and this central server. A native App for Smartphones allows EV users to search for free CSs, book CSs, set the preferences for the PEV recharging process (as mainly functionality), read the real time level of battery, the estimated autonomy of PEV, etc. The control center provides commands and load schedules to the CSs, which therefore only play the executive role in the whole energy management process, by adapting the charging load according to the instructions received.

The reservation process is carried out using the native smartphone App. The PEV user selects the most suitable CS and the most optimal FEV route is calculated and displayed on a map, also indicating the predicted arrival time, the route distance and the forecasted arrival state of charge. After having selected a charging station, the user introduces, via a web form interface, the needed parameters to allow the charging load control algorithm to compute the optimal charging profile for the charging sessions. These parameters include: (i) the desired final level of charge, (ii) the time when the charging should start, (iii) the time when the charging process has to be terminated and, (iv) the chosen socket/plug among the available ones (e.g. SCAME or Mennekes). Together with these user parameters, several selected FEV model characteristics are obtained from the SMARTV2G system information back-end: (i) maximum battery capacity, (ii) maximum and minimum accepted power, (iii) arrival status of charge and, (iv) number of phases are also stored within the SMARTV2G reservation management database in order to allow a correct configuration of the future charging session. Once the user has defined these parameters, the SMARTV2G system receives this info and connects to the selected CS to send this information. The authentication process verifies the identity of each user through an RFID card. The user code is automatically checked through the control center to see if the user is authorized to perform the charging of the vehicle and if a valid reservation for the user exists. In positive case, the control center creates a new charging session within the system database specifying the current start time (i.e. the first available timeslot for the load control algorithm to optimize) and the charging preferences (previously defined by the user during the reservation process described above).

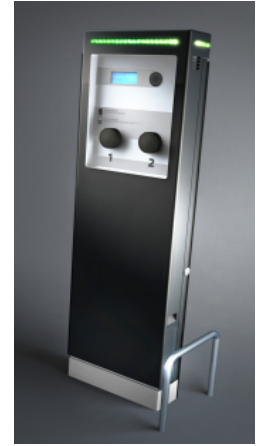

Fig. 3. Commercial image of Charging Station (2 x Type 2 socket-outlets)

Periodically (every 5 minutes), the charging station reports to the control center several information about the current charging sessions (e.g. applied power based on charging load set point calculated by the charging load control algorithm), which are used as feedback by the charging control algorithm to adjust the charging load set point if needed. When the charging session is terminated by the user (via his/her RFID card or by unplugging the vehicle), the charging station sends the last update to the control center (which terminates the charging session) and informs the user about the consumed energy in current charging session in kWh.

## IV. SMART CHARGING STATION

Smart charging station, used in the SMARTV2G field test, is developed and manufactured by ETREL (Fig. 3; see [11]). It is intended for public and semi-public use and enables simultaneous charging of two PEVs with maximum charging current of 3 x 32 A (22 kW) per socket-outlet. Charging station may be equipped with any combination of standard household socket-outlets (mode 2 charging according to IEC 61851 standard [10]) or Type 2 socket-outlets (in accordance with IEC 62196 standard) which enable mode 3 charging [10].

The main functional characteristics of CS used in SMARTV2G field tests are:

- Two Type 2 socket-outlets for mode 3 charging with maximum charging current of 3 x 32 A per socket-outlet;
- Utility feeder equipment (main disconnector/breaker, overcurrent and overvoltage protection, certified meters) embedded in the CS;
- Communication with PEV (PWM signal on control pilot wire) according to IEC 61851 standard [10] used also for PEV charging load control;
- PEV user identification with RFID cards that are swiped over the LED illuminated login area;
- LCD display for interaction with PEV user;

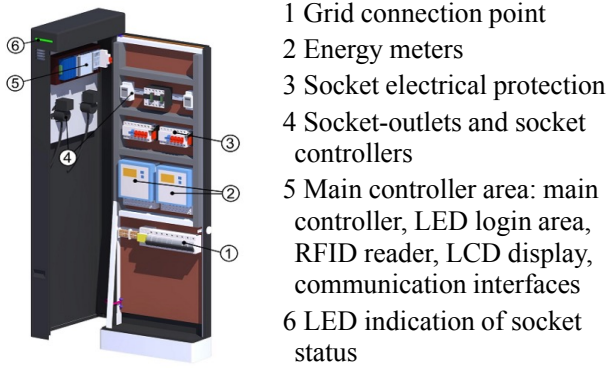

1 Grid connection point
2 Energy meters
3 Socket electrical protection
4 Socket-outlets and socket controllers
5 Main controller area: main controller, LED login area, RFID reader, LCD display, communication interfaces
6 LED indication of socket status

Fig. 4. Arrangement of equipment inside the charging station

- LED indicators to display the availability of CS;
- Energy meters for each socket acquire data on active (kWh) and reactive (kvarh) energy flow in both directions, output voltage (V) and current (A), and active power (kW);
- Bidirectional communication with SMARTV2G Control Centre for charging authorization, reporting and operation control and monitoring via web services using LAN/WAN or GSM network.

Arrangement of equipment inside the CS is presented in Fig. 4. The main control functions are executed by the CS main controller (Fig. 4, pos. 5) which communicates with all devices installed in the CS, controls their operation and communicates with SMARTV2G control center. The CS's main controller internally stores the charging load profile and appurtenant time intervals. During charging, each power set point to be applied at a defined time is converted, with consideration of number of charging phases, to current (A) integer value and communicated to socket controller. The socket controller further converts the set point current value to corresponding duty cycle of PWM signal on the control pilot and thus limits (controls) the PEV battery charger's load.

During charging, the CS sends the status report to the control center every 5 minutes. This information enables the control center to monitor the execution of charging load control and to take appropriate measures (modification of charging load profile) if PEV charging doesn't follow the load set points or if external conditions (new reference load profile given by DSO, new PEV(s) asking to charge) require the charging load reprofiling. At the end of charging session (repeated identification of PEV user and plugging out of the cable) the CS sends to control center a charging session report which contains detailed data about charging.

## V. LOAD CONTROL METHODOLOGY

As regards the demand side management algorithm, the charging sessions are managed by the LAC module according to a model predictive control principle [7]. Model predictive control is a popular optimization-based technique that allows to easily translate requirements and technical specifications of the problem into mathematical formulation. The reader is referred to [13]–[17] for a review of other possible approaches to smart charging.

The LAC module is in charge of updating the ongoing charging sessions (i.e. the load profiles to be applied by the EVSEs) each time a relevant event is notified to the EVSE operator, such as: new user requests, updates of the user preferences by one or more PEV users, price signals, volume signals by the DSO, etc. That is done in order to ensure that PEV charging load profiles are kept aligned with respect to the current state of the network and the current drivers' requirements. In the developed model predictive control framework, that is done by building, at each time of update, a new optimization problem based on the latest available information (i.e. number and type of charging sessions to be controlled, feedback from the CSs and the PEVs in the form, respectively, of metering data and current SoC values (if available), tariff updates, volume signals from the DSO, etc.).

### A. Mathematical Formulation

The mathematical formulation of the problem, which is the subject of [7], is briefly recalled in the following for the sake of completeness. The optimization problem is based on mixed-integer linear programming. The objective function is given by the cumulative costs for charging $J_{cost}$, plus the weighted distance $J_{reg}$ of the total PEV charging load from the DSO-defined reference load profile

$$J = J_{cost} + \mu J_{reg} \quad (1)$$

The cumulative costs for charging are given by

$$J_{cost} = \sum_{m \in M} \sum_{k=I}^{E_m-1} \Delta P_m T C[k] U_m[k] \quad (2)$$

where $k$ is the temporal index, $M$ is the set of controlled charging sessions, $I$ the time when the LAC is triggered for computation of the load profiles, $E_m$ the departure time of the $m$th PEV, $\Delta P_m$ the maximum power for the $m$th charging session, $T$ the discretization time, $C[k]$ the tariff and $U_m[k]$ the control signal to be applied by the EVSE. In compliance with standard IEC61851, $U_m[k]$ is a semi-continuous variable ($U_m[k] \in \{0\} \cup [\alpha_m, 1]$, with $0 < \alpha_m < 1$).

The load reference tracking term $J_{reg}$ is given by

$$J_{reg} = \left\| \Lambda (P - P_{ref}) \right\|_\infty \quad (3)$$

where $\Lambda$ is a diagonal matrix of weights, $P$ is the vector of aggregated load (whose $k$th component is the controlled load at time $k$, i.e. $P[k] = \sum_{m \in M_k} \Delta P_m U_m[k]$, where $M_k \subseteq M$ is the set of charging sessions controlled at time $k$) and $P_{ref}$ a vector whose $k$th component is the DSO defined load reference at time $k$. The operator $\|\cdot\|_\infty : R^n \to R$ denotes the the $l_\infty$- norm defined on the real space.

Several constraints are put in order to achieve a technically feasible solution of the optimization problem, in terms of real load profiles to be applied by the EVSEs. A constraint regards overload management, meaning that the cumulative controlled electromobility load cannot exceed a threshold $P^*$ set by the DSO for safe operation of the network

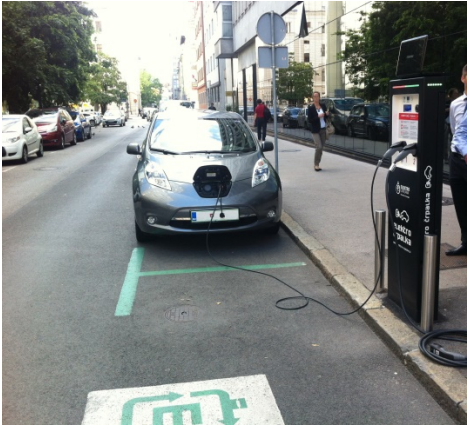

Fig. 5. Field test site with EVSE and PEV

$$\sum_{m \in M_k} \Delta P_m U_m[k] \leq P^*[k] \quad \forall k \in [I, E-1] \quad (4)$$

where $E = max\{E_m : m \in M\}$ is the last time of control problem definition. A second set of constraints takes into account the SoC limitation given by the battery capacity

$$X_m^{min} \leq x_m[k] \leq X_m^{max} \quad \forall m \in M \quad \forall k \in [I, E] \quad (5)$$

A third set of constraints makes sure that the user preferences are satisfied, namely, that the final state of charge is equal or greater than the one required by the PEV user

$$X_m^{ref} \leq x_m[E_m] \leq X_m^{max} \quad \forall m \in M \quad (6)$$

A simple model with constant losses $\xi_m$ has been used to predict the future SoC of controlled PEVs (which, notice, is needed in order to be able to write the above constraints)

$$\begin{cases} x_m[k+1] = x_m[k] + \Delta P_m T(1 - \xi_m) U_m[k] \\ x_m[I] = X_m^0 \end{cases} \quad (7)$$

Finally, a set of constraints is put in order to keep under control the costs associated to the single charging sessions (see [7] for a detailed explanation)

$$c_m[I] + \sum_{k=I}^{E_m-1} \Delta P_m T C[k] U_m[k] \leq (1+\varepsilon) c_m^* \quad \forall m \in M \quad (8)$$

The reader is referred to [7] for a complete explanation of the mathematical formulation of the control problem.

## VI. FIELD TEST RESULTS

Two kinds of tests have been performed in Ljubljana, Slovenia, at the premises of the Slovenian DSO Elektro Ljubljana and using the Elektro Ljubljana operated CSs:

- Test I: Application of a static charging load set point, computed after the charging request and never updated;
- Test II: Charging load set points periodically updated. The charging load set point is first computed after the charging request and then periodically updated until the end of the charging session. The charging load set point is updated to better track the power reference and to properly react to signals coming from the DSO and asking for a modification of the reference/maximum allowed load profile.

### A. Field Test Setup

The tests were performed on one of the Elektro Ljubljana CSs installed and operative on the field for public recharging, and by recharging a Nissan Leaf PEV (see Fig. 5 and TABLE I). The considered CS is supplied with the 3 phase connection, 3 x 63 A (main supply), voltage 230/400 V nominal power 43.65 kW, power factor 1.00, number of poles 3P+N, grounding TN-S. The main supply current is divided inside the EVSE for two sockets, so that inside the station two sockets are protected by 3 x 32 A protection. The consumption of the CS is measured with the installed smart meters for each socket separately. So two meters are installed. Such approach is specific but more economical because the smart meters cover the data collection for the DSO and also for the CS operator. So the DSO, as responsible for the consumption data collection, delivers these data to all actors needed (e.g. supplier, service provider). Modern smart meters enable also voltage level read out. When the procedure of read out is done, the data voltage is read out also.

### B. Test I: Application of Static Charging Load Set Points

This kind of test aimed at demonstrating on the field the correct functioning of the whole SMARTV2G PEV charging control system, from the correct handling of users' charging requests, up to computation of load profiles, correct application of profiles by the CS and finally the correct termination of the charging session. These tests have been performed with a single PEV and by letting the LAC compute the load profile only at the beginning of the charging session (i.e. without updating it periodically).

These kinds of tests have highlighted, first of all, that the knowledge of good estimates of the efficiency factor $\xi_m$ and, in particular, of the actual PEV battery capacity are

TABLE I. PEV PARAMETERS

| Nominal Battery Capacity [kWh] | 22 |
|---|---|
| Maximum charging power [kW] | 3.68[a] |
| Type of charging (single phase or three phase) | Single phase |

a. 1 * 16 A * 0.23 kV = 3.68 kW

TABLE II. TEST I

| Session ID | Initial SoC→Desired final SoC [%] | Measured final SoC [%] | Energy measured at CS [kWh] |
|---|---|---|---|
| 180[a] | 10→30 | 38 | 5.54 |
| 182[b] | 38→48 | 49 | 2.60 |

a. Load profile computed considering for the battery capacity the nominal value of 22 kWh
b. Load profile computed considering for the battery capacity the nominal value of 19 kWh

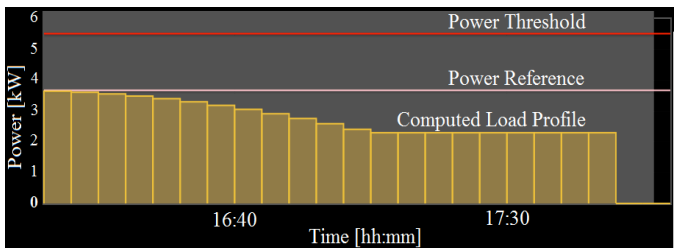

Fig. 6. Charging load set point as computed by the DSM algorithm

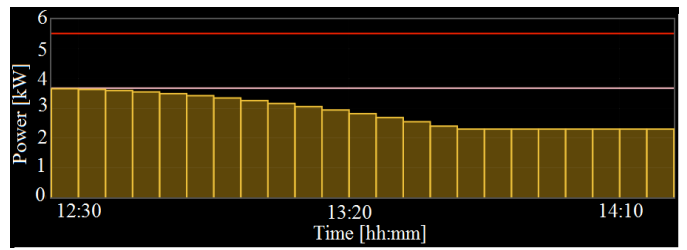

Fig. 8. First computed charging load set point

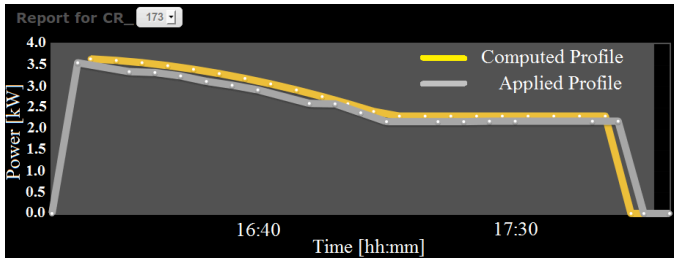

Fig. 7. Comparison between the computed load profile and the applied one

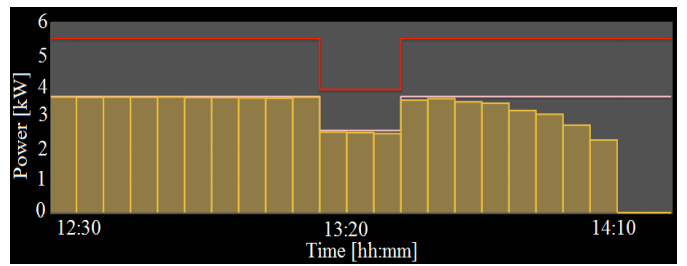

Fig. 9. Last computed charging load set point

fundamental to satisfy the requirement on user preferences (in particular, the requirement on the energy to recharge). That is because the demand side management algorithm expresses the SoC in kWh values, while the user preferences and the PEV dashboard readings express the SoC value in %, so that a conversion is needed, which in turn requires the knowledge of the actual SoC. TABLE II reports two type I tests. From the table it is seen that considering the nominal value of the battery capacity (first row) leads to overcharge the PEV: it has to be taken into account that battery degradation progressively reduces the actual battery capacity. From TABLE II, it can be seen that the results improve when a value of 19 kWh is considered for the actual battery capacity (notice that a first estimate for the actual battery capacity can be calculated by simply comparing the metering data – energy – from the CS and the SoC readings from the PEV).

As regards the requirement on the tracking of the reference load profile, Fig. 6 reports the relevant power signals derived from one of the performed tests (the power threshold, the power reference and the computed load profile, which is the load profile set point sent to the CS). It is seen that load reference tracking is accurate only close to the initial time (which is expected in these type I tests [7], since the load profile is never updated after it is first computed). Figure 7 compares the PEV charging load set point (computed profile) with the actual power applied by the charging station. The profile in grey is the one actually applied by the EVSE. Each sample (i.e. the white dots in the figure) of the applied profile represents the real power as measured by the CS's meter and sent to the SMARTV2G control centre each five minutes. From the figure it can be noticed that the applied power follows the calculated one with good accuracy. In particular, the average distance between the load profile set point and the actuated one is 0.114 kW, the mean squared error is 0.00372. The maximum distance between the load profile set point and the actuated one is 0.212 kW. The mismatch between the two profiles could depend on:

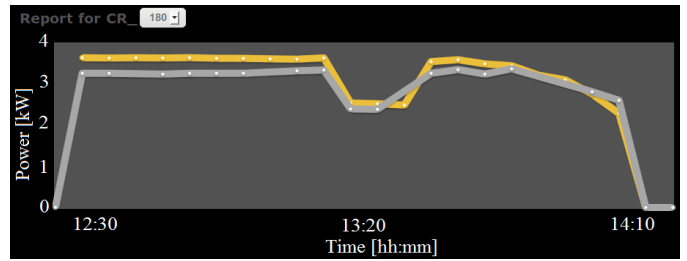

Fig. 10. Last computed charging load set point (in yellow) and applied load profile (in grey)

- Rounding of the power set point value to integer phase current values;
- Disturbances, like the ones influencing the voltage level, which is assumed constant by the DSM algorithm;
- Nonlinearities characterizing the battery and the charger.

### C. Test II: Charging Load Set Points Periodically Updated

In the second kind of field tests the control system periodically (each 5 minutes) updates the charging load set points, by iterating the EV charging control algorithm described in Section V. This approach allows achieving accurate reference tracking and to react to events like new charging requests and variations of the reference profile set by the DSO (i.e. DSM signals). As an example, in Fig. 8 and Fig. 9 are reported, respectively, the first and last charging load set points computed by the control system to serve a user charge request in such a "time-driven" mode. It can be seen that the load set point changes in order to best adapt to the load reference which, also, has been changed during the execution of the test, to simulate a volume reduction DSM signal sent by the DSO to limit the power available for charging. Figure 10 shows the final evolution of the charging load set point to the CS and the charging power actually applied, as measured by the CS meter. It is relevant to notice that some reports to the CS are lost due to GSM connectivity issues.

VII. CONCLUSIONS

This paper has discussed the real implementation of advanced electric vehicle recharging control functionalities supporting demand side management strategies. Field test results have been presented and discussed, demonstrating the feasibility and good performances of the proposed PEV charging management system. It has been demonstrated the ability of controlling the single load profiles at EVSE level according to requirements related to grid operation constraints, costs minimization and user preferences satisfaction. Also, the aggregated load from PEVs can be shaped in order to mitigate the impact on the grid and providing balancing power. Important remarks in particular relates to the necessity of accurate feedback from the PEVs regarding the actual capacity of the batteries and the real-time value of the SoC. The lack of such information precludes the possibility of accurately satisfying the user requirement on the final desired SoC.

Future works will regard the test of the proposed PEV control system for controlling simultaneous charging processes at public CSs, the test in a residential scenario, according to the PEV integration scheme proposed in [18], and the extension of the proposed control strategy for balancing renewable energy sources, especially dealing with photovoltaic plants and wind-turbine driven generators (the latter by proper interfacing with the control system developed in [19]). Novel Future Internet paradigms [20] and innovations in the ICT/telecommunication sectors [21] are also expected to bring benefits to the whole monitoring and control chain discussed in the paper and will be considered for future works.

ACKNOWLEDGMENT

The authors gratefully acknowledge Sixto Santonja Hernández and Caterina Tormo Domènech, ITE, Spain, for their precious coordination work in SMARTV2G project; Francesco Delli Priscoli and Francisco Facchinei, "Sapienza" University of Rome, Italy, for their valuable suggestions on the reference scenario and control system design and, finally, Zsolt Krémer and Andreas Varesi, Technomar GMBH, Germany, for their support in dissemination.

REFERENCES

[1] EU Commission. Communication of European commission to European parliament: an energy policy for Europe; 2007.

[2] Falvo, M.C.; Graditi, G.; Siano, P., "Electric Vehicles integration in demand response programs," Power Electronics, Electrical Drives, Automation and Motion (SPEEDAM), 2014 International Symposium on, pp. 548-553, 18-20 June 2014. DOI:10.1109/SPEEDAM.2014.6872126.

[3] M. C. Falvo, R. Lamedica, R. Bartoni, and G. Maranzano, Energy management in metro-transit systems: An innovative proposal toward an integrated and sustainable urban mobility system including plug-in electric vehicles, Electric Power Systems Research, Volume 81, Issue 12, December 2011, Pages 2127-2138, ISSN 0378-7796, DOI:10.1016/j.epsr.2011.08.004.

[4] SMARTV2G Consortium, 2013. Project website. URL ⟨www.SMARTV2G.eu/⟩.

[5] Di Giorgio, A.; Liberati, F.; Canale, S., "Optimal electric vehicles to grid power control for active demand services in distribution grids," Control & Automation (MED), 2012 20th Mediterranean Conference on, pp. 1309-1315, 3-6 July 2012. DOI: 10.1109/MED.2012.6265820.

[6] Di Giorgio, A.; Liberati, F.; Canale, S., "IEC 61851 compliant electric vehicle charging control in smartgrids," Control & Automation (MED), 2013 21st Mediterranean Conference on, pp. 1329-1335, 25-28 June 2013. DOI: 10.1109/MED.2013.6608892.

[7] Di Giorgio, A.; Liberati, F.; Canale, S., Electric vehicles charging control in a smart grid: A model predictive control approach, Control Engineering Practice, Volume 22, January 2014, Pages 147-162, ISSN 0967-0661. DOI: 10.1016/j.conengprac.2013.10.005.

[8] Di Giorgio, A.; Liberati, F., "Model Predictive Control-based approach for Electric Vehicles charging: power tracking, renewable energy sources integration and driver preferences satisfaction", in Plug-In Electric Vehicles in Smart Grid: Charging Strategies, Springer, 2014.

[9] Di Giorgio, A.; Liberati, F.; Pietrabissa, A., "On-board stochastic control of Electric Vehicle recharging," Decision and Control (CDC), 2013 IEEE 52nd Annual Conference on, pp.5710,5715, 10-13 Dec. 2013. DOI: 10.1109/CDC.2013.6760789.

[10] ISO/TC/SC 69 Electric road vehicles and electric industrial trucks. IEC 61851-1ed2.0 Electric vehicle conductive charging system – Part 1: General requirements; 2010.

[11] www.etrel.si; http://evcharging.etrel.com/

[12] IEC 62196-1 Plugs, socket-outlets, vehicle connectors and vehicle inlets – Conductive charging of electric vehicles – Part 2: Dimensional compatibility and interchangeability requirements for a.c. pin and contact-tube accessories.

[13] J. García-Villalobos, I. Zamora, J.I. San Martín, F.J. Asensio, V. Aperribay, Plug-in electric vehicles in electric distribution networks: A review of smart charging approaches, Renewable and Sustainable Energy Reviews, Volume 38, October 2014, Pages 717-731, ISSN 1364-0321. DOI: 10.1016/j.rser.2014.07.040.

[14] Sortomme, E.; El-Sharkawi, M.A., "Optimal Charging Strategies for Unidirectional Vehicle-to-Grid," Smart Grid, IEEE Transactions on, vol.2, no.1, pp.131,138, March 2011. DOI: 10.1109/TSG.2010.2090910.

[15] Deilami, S.; Masoum, A.S.; Moses, P.S.; Masoum, M.A.S., "Real-Time Coordination of Plug-In Electric Vehicle Charging in Smart Grids to Minimize Power Losses and Improve Voltage Profile," Smart Grid, IEEE Transactions on, vol.2, no.3, pp.456,467, Sept. 2011. DOI: 10.1109/TSG.2011.2159816.

[16] Masoum, A.S.; Deilami, S.; Moses, P.S.; Masoum, M.A.S.; Abu-Siada, A., "Smart load management of plug-in electric vehicles in distribution and residential networks with charging stations for peak shaving and loss minimisation considering voltage regulation," Generation, Transmission & Distribution, IET, vol.5, no.8, pp.877,888, August 2011. DOI: 10.1049/iet-gtd.2010.0574

[17] Masoum, Amir S.; Deilami, Sara; Masoum, Mohammad A.S.; Abu-Siada, Ahmed; Islam, Syed, "Online coordination of plug-in electric vehicle charging in smart grid with distributed wind power generation systems," PES General Meeting | Conference & Exposition, 2014 IEEE, vol., no., pp.1,5, 27-31 July 2014. DOI: 10.1109/PESGM.2014.6939133

[18] Di Giorgio, A.; Liberati, F., Near real time load shifting control for residential electricity prosumers under designed and market indexed pricing models, Applied Energy, Volume 128, 1 September 2014, Pages 119-132, DOI: 10.1016/j.apenergy.2014.04.032.

[19] Di Giorgio, A.; Pimpinella, L.; Mercurio, A., "A feedback linearization based Wind turbine control system for ancillary services and standard steady state operation," Control & Automation (MED), 2010 18th Mediterranean Conference on , vol., no., pp.1585,1590, 23-25 June 2010. DOI: 10.1109/MED.2010.5547821.

[20] M. Castrucci, F. Delli Priscoli, A. Pietrabissa, V. Suraci, "A Cognitive Future Internet Architecture", in "The Future Internet" - second edition, Springer, 2011, ISBN 978-3-642-20897-3, pp. 91-102.

[21] Canale, S.; Di Giorgio, A.; Lanna, A.; Mercurio, A.; Panfili, M.; Pietrabissa, A., "Optimal Planning and Routing in Medium Voltage PowerLine Communications Networks," Smart Grid, IEEE Transactions on, vol.4, no.2, pp.711-719, June 2013. DOI: 10.1109/TSG.2012.2212